# Virtual Power Plant (VPP) Architecture Tradeoffs

Deepak Aswani
Sacramento Municipal Utility District
Sacramento, CA, USA
Deepak.Aswani@smud.org

*Abstract*—With a broad market of Distributed Energy Resource (DER) aggregated control solutions, technology investment and implementation decisions for electric utilities and grid operators need to consider complexity, cost, and performance. This paper compares the performance tradeoffs of two Virtual Power Plant (VPP) architectures for DER aggregated control.

## I. Background

In 2021, the board of the Sacramento Municipal Utility District (SMUD) adopted an ambitious goal to decarbonize its energy supply by 2030. With a shift in the energy supply portfolio to variable renewable generation, a need for substantial amount of energy storage and load flexibility was anticipated to balance energy supply with demand.

Customer Distributed Energy Resources (DER) in load flexibility pilots have demonstrated that uncertainty and variability in performance is significant which has led to innovation in Virtual Power Plants (VPP) that attempt to firm load response. Control of DER in the form of VPPs may be a more cost-effective solution than utility battery storage at low utilization rates because it shares customer investments. The market for VPP products and DER aggregation is rapidly expanding with many alternatives to consider.

One specific tradeoff for VPP solutions is whether to pursue layered aggregation with many DER-specialized VPP operators or a more centralized aggregation involving few VPP operators as depicted in Figure 1. Products and services are available in the market today for both architectures. Comparison of these two approaches will be made where each aggregation is subject to the same optimization objective and constraints using data from real world experiments.

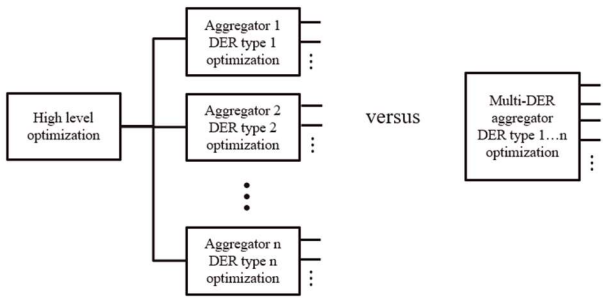

Figure 1: Two architectures compared

The centralized approach can simplify Information Technology (IT) investments for back-office operational integration, but the performance tradeoffs of both architectures are largely unknown. This paper uses experimental data from pilots conducted at SMUD with DER aggregators to model the difference in performance under these two architectures in terms of flexible load/generation power capability and statistical certainty over time.

## II. Experiment Data

Direct experimental results from residential customer pilots were used for this analysis unless otherwise specified:

- Smart thermostats – Demand response for Nest and Ecofactor smart thermostats [1].
- Electric Vehicles – ClipperCreek smart Level 2 chargers (EVSE) [2]. The model of [3] was fit to this data to simulate additional experiments.
- Water heaters – GE and Rheem controllable water heaters managed by aggregator Virtual Peaker [4].
- Batteries – Controlled SolarEdge residential battery energy storage systems capable of grid export [5].

## III. DER Modeling Framework

Load flexibility for an aggregation of DER can be modeled as a time-indexed convex space bound by periodic load shapes that are a function of customer/utility control inputs, weather conditions, and starting state.

The dynamic response of a group of DERs can be generalized using a discrete-time non-linear state-space representation as described as the foundational plant model for any non-linear controls problem in [6]:

$$\mathrm{x}[k+1] = f(k, \mathrm{x}[k], \mathrm{u}[k], \mathrm{w}[k]) \quad (1)$$

$$y_{Power}[k] = g(k, \mathrm{x}[k], \mathrm{u}[k], \mathrm{w}[k]) \quad (2)$$

where both $f(\cdot)$ and $g(\cdot)$ are non-linear functions and

$\mathrm{x}[k]$ = time-indexed state-space column matrix of all DER,

$\mathrm{u}[k]$ = time-indexed column matrix of control commands for the DER,

$\mathrm{w}[k]$ = time-indexed column matrix of uncontrollable inputs (i.e., ambient temperature, solar irradiance, etc.), and

$y_{Power}[k]$ = total power of all DER at time index $k$ relative to baseline.

Given diurnal patterns of customer use of DER and weather, one may assume that the state of all DER in aggregate approximately resets daily (i.e., at midnight, 4 AM, etc.). Examples include the fact that Electric Vehicles are expected to be fully charged by the daily commute departure, most battery systems are expected to be sufficiently depleted to allow for day-time solar recharge but not overly depleted to compromise backup power reserve, etc. However, there is evidence that two days are necessary to capture periodic patterns due to prior day state dependence [5]. For this smallest periodic time interval, the reset state of DER is assumed to be approximately constant.

$$X[k_{\text{reset hour(day+1)}}] \approx x[k_{\text{reset hour(day)}}] \approx C \quad (3)$$

Suppose there are $N_{Period}$ steps to the periodic sequence. For ease of notation, it is assumed that $k$ is offset such that every multiple of $N_{Period}$ is a reset time (i.e. $N_{Period}, 2N_{Period}, 3N_{Period}, \ldots$). For the period following period $i$, iterative substitution of (1), (2), and (3) can be simplified to provide an important insight – the total DER power in aggregate can be simplified as a new function of daily sequences that is dependent on starting state as shown in (4).

$$\vec{y}_{Power}(\{D\})[i] = h(i, C, \vec{u}[i], \vec{w}[i]) \quad (4)$$

This equation (4) concerns period $i$ for a presumed set $D$ of DER, where the daily sequence inputs and outputs beginning at the time of the daily reset are defined as

$$\vec{u}[i] = [u[iN_{Period}], \ldots, u[(i+1)N_{Period} - 1]] \quad (5)$$

$$\vec{w}[i] = [w[iN_{Period}], \ldots, w[(i+1)N_{Period} - 1]] \quad (6)$$

$$\vec{y}_{Power}(\{D\})[i] = [y_{Power}[iN_{Period}], \ldots, y_{Power}[(i+1)N_{Period} - 1]] \quad (7)$$

For addition of DER of the same mix (device, manufacturer, connectivity provider, etc.), a multiplicative scaling of the total aggregate power is expected. For example, if the existing mix of DER is doubled, then double the mean response would be expected in aggregate if the inputs $\vec{u}[i]$ and $\vec{w}[i]$ are the same therefore having the same mean impact in aggregate for a large population of DER. Randomly sampled fraction $\alpha$ of $\{D\}$ or

$$\{D\} \xleftarrow{R} \{D_\alpha\} \quad (8)$$

will result in a multiplicative scaling of (4) that can be extended as

$$\vec{y}_{Power}(\{D_\alpha\})[i] = \alpha h(i, C, \vec{u}_\alpha[i], \vec{w}[i]) \quad (9)$$

Partitioning the DER into $Q$ groups with random assignment according to (9) results in

$$\vec{y}_{Total} = \sum_{j=1}^{Q} \alpha_j \vec{y}_{Power,j} \quad (10)$$

subject to the constraints

$$\sum_{j=1}^{Q} \alpha_j = 1 \text{ and } \alpha_j \geq 0 \quad (11)$$

which means that $\{\vec{y}_{Total}\}$ is a convex set.

IV. EXPERIMENT DESIGN

To define the DER convex set, various permutations of customer/utility control inputs, weather conditions, and starting state must be considered as in Figure 2.

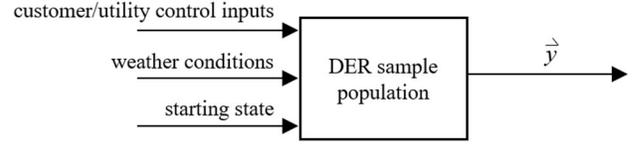

Figure 2: DER black box for experimental data collection

The volume of data collected should include sufficient permutations of these inputs. The customer/utility control inputs should be orthogonal to get sufficient variability but also be representative of sequences that would be considered reasonable to the utility and customer. Supposing that output measurements statistically summarized by unique combinations of input result in a set of vector pairs of mean load shape and associated variance over the periodic interval, nominally 48 hours for a single DER:

$$\left\{ \begin{bmatrix} \vec{y}_\mu \\ \vec{y}_{\sigma^2} \end{bmatrix} \right\} \quad (12)$$

From this set of vector pairs, numerical algorithms can be used to identify the convex hull for $\vec{y}_\mu$ where

$$\{\vec{y}_{\mu Hull}\} = \text{Hull}(\vec{y}_\mu) \in \{\vec{y}_\mu\} \quad (13)$$

of size $Q_{Hull}$. The set $\{\vec{y}_{\mu Hull}\}$ yields in a reduced set of vector pairs based on the original vector pair association (12) defined as

$$\left\{ \begin{bmatrix} \vec{y}_{\mu Hull} \\ \vec{y}_{\sigma^2_*} \end{bmatrix} \right\} \quad (14)$$

It is important to note that only $\{\vec{y}_{\mu Hull}\}$ is a convex subset according to Section III. The set $\{\vec{y}_{\sigma^2_*}\}$ does not exhibit convex properties and needs special consideration as explained in the following performance framework section.

DER data was collected in experiments identified in Section II as vector pair sets defined in (12) that represent and the real world variability and uncertainty in performance. Given the limited time frame of testing, vector pair data were collected for 24-hour periodic intervals rather than 48-hour periodic intervals to expand the sample size.

V. PERFORMANCE FRAMEWORK FOR DER CONTROL

The tradeoff of DER flexible load/generation power capability and reliability can be balanced via co-optimization.

If a portfolio of DER is intended to operate as a Virtual Power Plant (VPP), flexible load/generation power capability and reliability can be described as load shift relative to baseline and statistical confidence of that load shift, respectively.

The technical potential relative to baseline for time index $k$ can be characterized as

$$\min\{\text{Bound} \cdot \text{mean}(\vec{y}_{Net}[k]) + \lambda \cdot \text{variance}(\vec{y}_{Net}[k])\} \quad (15)$$

where $\lambda$ is a scalar weighting factor that adjusts the tradeoff between mean and variance and Bound is defined as -1 and 1 to characterize both upper and lower bound technical potential respectively. Suppose that $N_{Total}$ describes the total number of similar DER type (i.e., EV, battery, water heater, thermostat, etc.), so that dispatch of a homogeneous population of DER is partitioned across the $D_{Hull}$ points of the convex hull

$$N_{Total} = \sum_{j=1}^{D_{Hull}} n_j \tag{16}$$

where noting that in (9), $i$ drops out as it reflects the periodic interval, $\vec{u}_\alpha$ is a function of $\alpha$ which is analogous to being indexed by $j$, from here on referenced as , $\vec{u}_j$ then

$$\text{mean}(\vec{y}_{Net}[k]) = \sum_{j=1}^{D_{Hull}} \left(n_j \cdot \vec{y}_{\mu Hull,j}(C, \overline{w})[k]\right) \tag{17}$$

Identical distribution may only be unique to each partition associated with a point on the convex hull per (14). Statistical independence is not guaranteed because of coupling to common conditions $\overline{w}$, starting state $C$, and control inputs $\vec{u}_j$. The variance of each partition can therefore be defined as

$$\text{variance}(\vec{y}_{Net,j}[k]) = \sum_{p=1}^{n_j} \vec{y}_{\sigma^2*,j}(C, \overline{w})[k] + \cdots \tag{18}$$

$$2 \sum_{1 \leq a < b \leq n_j} \text{covariance}\left(\vec{y}_{j,a}(C, \overline{w})[k], \vec{y}_{j,b}(C, \overline{w})[k]\right)$$

Using the Pearson's correlation coefficient $\rho$, assume

$$\text{covariance}\left(\vec{y}_{j,a}(C, \overline{w})[k], \vec{y}_{j,b}(C, \overline{w})[k]\right) \tag{19}$$
$$= \rho_j(C, \overline{w})[k] \cdot \vec{y}_{\sigma^2*,j}(C, \overline{w})[k]$$

which after substitution in (18) and simplification yields

$$\text{variance}(\vec{y}_{Net,j}[k]) = \tag{20}$$
$$\left(n_j^2 \cdot \rho_j(C, \overline{w})[k] + n_j \cdot (1 - \rho_j(C, \overline{w})[k])\right)$$
$$\cdot \vec{y}_{\sigma^2*,j}(C, \overline{w})[k]$$

Across different $j$, since the control input $\vec{u}_j$ is unique, statistical independence is effectively assumed so that

$$\text{variance}(\vec{y}_{Net}[k]) = \tag{21}$$
$$\sum_{j=1}^{D_{Hull}} \left(\left(n_j^2 \cdot \rho_j(C, \overline{w})[k] + n_j \cdot (1 - \rho_j(C, \overline{w})[k])\right) \cdot \vec{y}_{\sigma^2*,j}(C, \overline{w})[k]\right)$$

Defining

$$L = \text{Bound} \cdot \left[\vec{y}_{\mu Hull,1}(C, \overline{w})[k], \ldots, \vec{y}_{\mu Hull,D_{Hull}}(C, \overline{w})[k]\right] \tag{22}$$
$$+ \lambda \cdot \left[(1 - \rho_1(C, \overline{w})[k]) \cdot \vec{y}_{\sigma^2*,1}(C, \overline{w})[k], \ldots, (1 - \rho_{D_{Hull}}(C, \overline{w})[k]) \cdot \vec{y}_{\sigma^2*,D_{Hull}}(C, \overline{w})[k]\right]$$

$$Q = \lambda \cdot \text{diag}\Big(\rho_1(C, \overline{w})[k] \tag{23}$$
$$\cdot \vec{y}_{\sigma^2*,1}(C, \overline{w})[k], \ldots, \rho_{D_{Hull}}(C, \overline{w})[k]$$
$$\cdot \vec{y}_{\sigma^2*,D_{Hull}}(C, \overline{w})[k]\Big)$$

$$n = [n_1, \ldots, n_{D_{Hull}}]^T \tag{24}$$

Then substitution of (17) and (21) into (15) and referencing (16) results in a quadratic optimization problem

$$\min_{n}\{L\,n + n^T Q\,n\} \text{ subject to} \tag{25}$$
$$1_{1,D_{Hull}}\,n = N_{Total} \text{ and} \tag{26}$$
$$n_j \geq 0 \ \forall\ j \in \{1, \ldots, D_{Hull}\} \tag{27}$$

where the optimization sizes each $n_j$ of $D_{Hull}$ partitions that is cumulatively $N_{Total}$ DER and each $n_j$ is uniquely associated with control sequence $\vec{u}_j$.

## VI. Extending the Performance Framework to Multiple DER Types

The mean load modification of aggregated DER can be represented as a convex set that is a Minkowski sum of subgroups of DER.

Section V presumes a homogeneous population of DER type. For a heterogeneous population of multiple types of DERs, each DER type can be described as having its own convex set $\{\vec{y}_{\mu Hull}(r)\}$ for each DER type $r$ of $m$ types. This is graphically depicted in Figure 3.

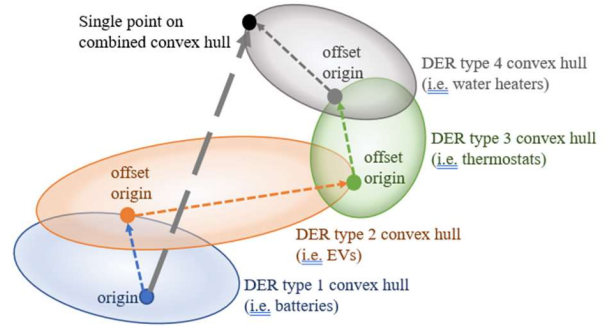

Figure 3: Example of a single point on the convex hull of multiple heterogenous DER

The Minkowski sum of the individual homogeneous hulls scaled to the volume of DER is

$$\text{Hull}\left(\sum_{r=1}^{m} n(r) \cdot \vec{y}_\mu(r)\right) = \sum_{r=1}^{m} \text{Hull}\left(n(r) \cdot \vec{y}_\mu(r)\right) \tag{28}$$
$$= \sum_{r=1}^{m} n(r) \cdot \{\vec{y}_{\mu Hull}(r)\}$$

This allows the model of Section V to extend to a heterogeneous population of DER. It is assumed that the performance of different DER types is statistically independent resulting in

$$L_{\text{multi}} = [L(1), \ldots, L(m)] \tag{29}$$

$$Q_{multi} = \text{diag}(Q(1), \ldots, Q(m)) \quad (30)$$

$$n_{multi} = \begin{bmatrix} n(1) \\ \vdots \\ n(m) \end{bmatrix} \quad (31)$$

With combined optimization

$$\min_{n_{multi}} \{L_{multi}\, n_{multi} + n_{multi}^T Q_{multi}\, n_{multi}\} \text{ subject to} \quad (32)$$

$$\begin{bmatrix} \mathbf{1}_{1,D_{Hull(1)}} & & 0 \\ & \ddots & \\ 0 & & \mathbf{1}_{1,D_{Hull(m)}} \end{bmatrix} n_{multi} = \begin{bmatrix} N_{Total}(1) \\ \vdots \\ N_{Total}(m) \end{bmatrix} \text{ and} \quad (33)$$

$$n_j \geq 0 \;\; \forall\; j \in \left\{1, \ldots, \sum_{r=1}^{m} D_{Hull(r)}\right\} \quad (34)$$

## VII. COMPARISON OF PERFORMANCE ENVELOPE BY ARCHITECTURE

Given the limited volume of data collection, it is not possible to estimate $\rho_j[k]$ from (19). However, analysis of load flexibility experiments in [7] found that for 95% of the experiments, the correlation of load shift across time was about 0.11 or less. Therefore, for simplicity a constant value of 0.1 is assumed for $\rho_j[k]$ with the following population of DER:

- 300 residential batteries
- 3,000 smart thermostats and controllable water heaters (each)
- 1,500 EVs workplace and residential charging (each)

The two high level architectures for VPP analyzed are from Figure 1. This is analogous to separate optimizations or multiple instances of the optimization problem described in (22)-(27) or a single instance of the optimization problem described in (29)-(34), respectively.

A centralized aggregation has greater flexible load/generation power capability than layered aggregation when tracking power while considering cost constraints.

When considering a VPP for displacing traditional generation, the optimization cost weighting for variance ($\lambda$) should balance flexible load/generation power capability and reliability. Dispatch cost is a consideration because operational decisions often find different generators and market transaction needs that affect the energy cost at the margin. Furthermore, with a likely model that DER are over-procured and not all DER are dispatched every day, the cost of VPP dispatch can vary. This changes the earlier equality constraints to inequality constraints. For layered aggregation, (26) becomes

$$\mathbf{1}_{1,D_{Hull(r)}}\, n(r) \leq N_{Total}(r) \text{ for } r \in \{1, \ldots, m\} \quad (35)$$

and for centralized aggregation, (33) becomes

$$\begin{bmatrix} \mathbf{1}_{1,D_{Hull(1)}} & & 0 \\ & \ddots & \\ 0 & & \mathbf{1}_{1,D_{Hull(m)}} \end{bmatrix} n_{multi} \leq \begin{bmatrix} N_{Total}(1) \\ \vdots \\ N_{Total}(m) \end{bmatrix} \quad (36)$$

Consideration of total daily budgeted cost (Budget) and daily unit dispatch cost ($\text{Cost}_{UD}$) introduces additional constraints. For layered aggregation this is

$$\text{Cost}_{UD}(r) \cdot \mathbf{1}_{1,D_{Hull(r)}}\, n(r) \leq \beta(r) \cdot \text{Budget} \quad (37)$$

where $\beta(r)$ is the fraction of budgeted cost per homogeneous DER type $r$ so that

$$\sum_{r=1}^{m} \beta(r) = 1 \quad (38)$$

For centralized aggregation the cost constraint is

$$\begin{bmatrix} \text{Cost}_{UD}(1) \cdot \mathbf{1}_{1,D_{Hull(1)}} & \cdots & \text{Cost}_{UD}(m) \cdot \mathbf{1}_{1,D_{Hull(m)}} \end{bmatrix} \leq \text{Budget} \quad (39)$$

When considering reliability, dispatch cost, and flexible load/generation power capability, layered and centralized aggregation perform differently. Assuming a dispatch cost between \$0.15-\$0.50 per device per event with a daily dispatch budget of \$500 and cost weighting for variance ($\lambda$) at 0.1, Figure 4 shows that centralized aggregation (solid line) outperforms layered aggregation (dotted line) as the VPP maximum ability to flex load power up and the maximum ability to flex generation power (absolute value) is greater for almost every hour, when subject to the same optimization objective and constraints.

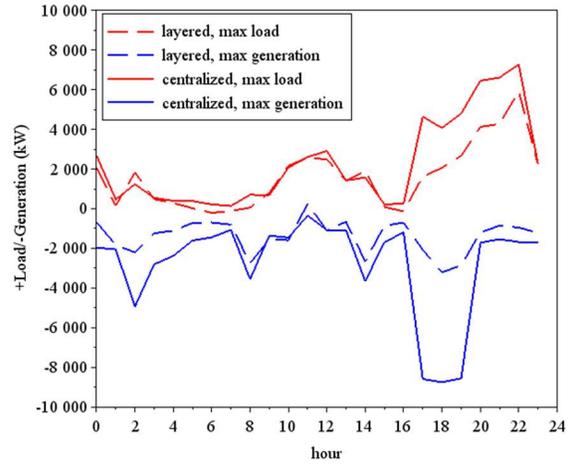

Figure 4: VPP performance envelope by aggregation architecture at 95% confidence level

## VIII. VPP IMPLEMENTAITON FOR POWER TRACKING AND ADAPTATIVE FEEDBACK

The optimal convex hull for centralized aggregation determined as in (29)-(34) can be summarized as

$$\Upsilon = \begin{bmatrix} \vec{y}_{Power,1} \\ \vdots \\ \vec{y}_{Power,Q} \end{bmatrix} \quad (40)$$

with associated control sequences $\{\vec{u}[i]\}$ and convex partition $\boldsymbol{\Lambda} = [\lambda_1, \ldots, \lambda_Q]$ Tracking $\vec{y}_{Total}$ with a day ahead power reference signal $\vec{y}_{Ref}$ can be achieved by minimizing mean squared error. The following minimization across $\boldsymbol{\Lambda}$ can be solved as a linear least squares problem:

$$\arg\min_{\boldsymbol{\Lambda}} \|\boldsymbol{\Lambda}\boldsymbol{\Upsilon} - \vec{y}_{Ref}\|_2^2 \qquad (41)$$
$$= \arg\min_{\boldsymbol{\Lambda}}(\boldsymbol{\Lambda}\boldsymbol{\Upsilon}\boldsymbol{\Upsilon}^T\boldsymbol{\Lambda}^T - 2\vec{y}_{Ref}\boldsymbol{\Upsilon}^T\boldsymbol{\Lambda}^T)$$

Complete visibility of measured or estimated power of individual DER allows for reconstructing the scaled $\vec{y}$ vectors based on known partitions and associated $\lambda$ values. These measurements are denoted as $\mathbf{Z}_{meas}$ shown in (42).

$$\mathbf{Z}_{meas} = \begin{bmatrix} \vec{y}_{scaled,q(1)} \\ \vdots \\ \vec{y}_{scaled,q(N)} \end{bmatrix} = \text{diag}(\boldsymbol{\Lambda}_{red})(\boldsymbol{\Upsilon}_{actual}) + \boldsymbol{\varepsilon} \qquad (42)$$

Directly solving for $\boldsymbol{\Upsilon}_{actual}$ is a challenge because $\text{diag}(\boldsymbol{\Lambda}_{red})$ can be ill-conditioned or close to singular. A low-pass discrete-time filter applied to both sides of (42) can reduce the ill-condition risk because $\boldsymbol{\Lambda}_{red}$ consists of only non-negative values per (11) and therefore the smallest values in the daily variation in $\text{diag}(\boldsymbol{\Lambda}_{red})$ will increase.

$$(\mathbf{Z}_{meas} * h_{filt})[i] = (\text{diag}(\boldsymbol{\Lambda}_{red}) * h_{filt})[i]\, \boldsymbol{\Upsilon}_{actual} \qquad (43)$$
$$+ (\boldsymbol{\varepsilon} * h_{filt})[i]$$

If $\boldsymbol{\varepsilon}$ is close to zero mean, then $(\boldsymbol{\varepsilon} * h_{filt})[i] \approx \mathbf{0}$. Defining $\mathbf{Z}_{meas}^{filt}[i] = (\mathbf{Z}_{meas} * h_{filt})[i]$ and $\text{diag}(\boldsymbol{\Lambda}_{red}^{filt}[i]) = (\text{diag}(\boldsymbol{\Lambda}_{red}) * h_{filt})[i]$ allows simplification of (43) to

$$\hat{\boldsymbol{\Upsilon}}_{actual} \approx \left(\text{diag}(\boldsymbol{\Lambda}_{red}^{filt}[i])\right)^{-1} \mathbf{Z}_{meas}^{filt}[i] \qquad (44)$$

A first order low-pass filter using a 1st order bilinear transform could be used resulting in

$$\mathbf{Z}_{meas}^{filt}[i] = \left(2\frac{\tau_{adapt}}{t_s} + 1\right)^{-1}\left(\left(2\frac{\tau_{adapt}}{t_s} - 1\right)\mathbf{Z}_{meas}^{filt}[i-1] + \mathbf{Z}_{meas}[i] + \mathbf{Z}_{meas}[i-1]\right) \qquad (45)$$

$$\text{diag}(\boldsymbol{\Lambda}_{red}^{filt}[i]) = \left(2\frac{\tau_{adapt}}{t_s} + 1\right)^{-1}\left(\left(2\frac{\tau_{adapt}}{t_s} - 1\right)\text{diag}(\boldsymbol{\Lambda}_{red}^{filt}[i-1]) + \text{diag}(\boldsymbol{\Lambda}_{red}[i]) + \text{diag}(\boldsymbol{\Lambda}_{red}[i-1])\right) \qquad (46)$$

where $t_s = 1$ the daily sampling period in days for $\mathbf{Z}_{meas}$ and $\tau_{adapt}$ = adaptation time constant in days.

Power tracking relative to reference is simulated along with a seven-day adaptation time constant over 60 days in Figure 5. The modeled feedback includes a 2% time-varying random error according to (42). Furthermore, at day 30 a significant change in DER model behavior is simulated to gauge the ability to adapt to disturbances and uncertainty.

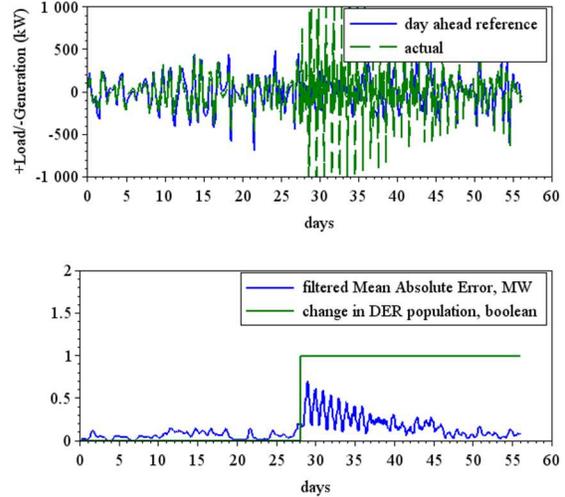

Figure 5: Power tracking and adaptive response to major disturbance/change in DER model behavior

## IX. CONCLUSION

Through the analysis of this paper using DER experiment data, it was determined that a centralized aggregation has greater flexible load/generation power capability than layered aggregation when tracking a power reference. Operation of such a VPP was simulated in response to disturbances. This DER architecture comparison for VPPs assumes a common optimization objective function and constraints. Such insight on performance can be helpful for utilities and grid operators to navigate tradeoffs for investments in DER control solutions and integration with back-office systems.

ACKNOWLEDGMENT

The author gratefully thanks the pilot data contributors L. Jimenez, A. Steeves, B. Harris, D. Huston, D. MacCurdy, N. Wong, B. Korven, O. Howlett, Y. Herrera, J. Duvall, J. Starrh, J. Frasher, D. Hinds, ClipperCreek, SolarEdge, and Virtual Peaker. The author also gratefully thanks the reviewers of this paper for improving the quality of content and presentation.